\documentclass[10pt]{article}
\usepackage{amssymb}
\usepackage{amsmath}
\usepackage{amsthm}
\theoremstyle{plain}

\newtheorem{theorem}{Theorem}[section]

\newtheorem{lemma}[theorem]{Lemma}
\newtheorem{corollary}[theorem]{Corollary}
\theoremstyle{definition}
\newtheorem{definition}[theorem]{Definition}

\theoremstyle{remark}
\newtheorem{remark}[theorem]{Remark}

\newcommand{\N}{\mathbb{N}}

\newcommand{\Z}{\mathbb{Z}}

\newcommand{\T}{\mathbb{T}}

\newcommand{\cM}{{\cal M}}

\newcommand{\osigma}{\overline{\sigma}}
\newcommand{\balpha}{\overline{\alpha}}

\newcommand{\oW}{\overline{W}}

\DeclareMathOperator{\id}{id}

\DeclareMathOperator{\oInn}{\overline{Inn}}
\DeclareMathOperator{\Aut}{Aut}
\DeclareMathOperator{\Ad}{Ad}
\DeclareMathOperator{\WInn}{WInn}

\newcommand{\bprf}{\noindent{\it Proof.\ }}

\newcommand{\eprf}{\hspace*{\fill} \rule{1.6mm}{3.2mm} \vspace{1.6mm}}
\newcommand{\benu}{\begin{enumerate}\renewcommand{\labelenumi}{{\rm (\roman{enumi})}}\renewcommand{\itemsep}{0pt}}
\newcommand{\eenu}{\end{enumerate}}
\begin{document}
\title{Discrete amenable group actions on von Neumann algebras and invariant nuclear C$^*$-subalgebras}
\author{Yasuhiko Sato \\}
\date{\small Department of Mathematics, Kyoto University\\
Sakyo-ku, Kyoto 606-8502, Japan \\
e-mail : ysato@math.kyoto-u.ac.jp}

\maketitle

\begin{abstract}   
Let $G$ be a countable discrete amenable group, $\cM$ a McDuff factor von Neumann algebra, and $A$ a separable nuclear weakly dense C$^*$-subalgebra of $\cM$. We show that if two centrally free actions of $G$ on $\cM$ differ up to approximately inner automorphisms then they are outer conjugate by an approximately inner automorphism, in the operator norm topology, which makes $A$ invariant. In addition, when $A$ is unital, simple, and with a unique tracial state and $\alpha$ is an automorphism of $A$ we also show that the aperiodicity of $\alpha$ on the von Neumann algebra is equivalent to the weak Rohlin property. 


\end{abstract}
\section{Introduction}
Since Connes' classification \cite{Con} of automorphisms of the injective type $I\hspace{-.1em}I_1$ factor von Neumann algebra, the classification of group actions on von Neumann algebras has been intensively studied. In particular, discrete amenable group actions on injective factors were completely classified by many hands. For finite group actions on the injective type $I\hspace{-.1em}I_1$ factor the classification was obtained by Jones \cite{Jon}, after that for discrete amenable group actions on type $I\hspace{-.1em}I$ factors this work was extended by Ocneanu \cite{Ocn}. For type $I\hspace{-.1em}I\hspace{-.1em}I_{\lambda}$ ($\lambda\neq 1$) factors the classification of discrete amenable group actions was obtained by Sutherland and Takesaki \cite{ST}, for type $I\hspace{-.1em}I\hspace{-.1em}I_ 1$ the classifications of finite group and abelian group actions were obtained by Kawahigashi, Sutherland, and Takesaki \cite{KwST}, and finally the classification was completed by Katayama, Sutherland, and Takesaki \cite{KST}. Recently, Masuda presented the unified proof for these results which is independent of types on factors \cite{Mas}, \cite{Mas2}.

In the proof shown by Masuda, he applied the Evans-Kishimoto intertwining argument for discrete amenable group actions on von Neumann algebras. The Evans-Kishimoto intertwining argument was first introduced in \cite{EK} to classify automorphisms of approximately finite dimensional C$^*$-algebras. Subsequent to this successful method, many classification results were obtained for amenable group actions on C$^*$-algebras \cite{Iz}. The purpose of the present paper is to reimport Masuda's Evans Kishimoto intertwining argument  from injective von Neumann algebras into nuclear C$^*$-algebras.

Our proof of the main result Theorem \ref{Thm1} is a combination of the intertwining argument and the amenability of nuclear C$^*$-algebras which was obtained by Haagerup in \cite{Ha}. In Section \ref{Sec2}, we shall obtain central sequences in the operator norm sense by the approximate diagonal of amenable C$^*$-algebras which was defined by Johnson \cite{Joh}. As an application we also show the second main result Theorem \ref{Thm2}. In Section \ref{Sec3}, applying Lemma \ref{Lem0} we shall follow the proof of Theorem 3 in \cite{Mas} and prove the second main theorem Theorem \ref{Thm1}.

Concluding this section, we prepare some notations. 
 For a $C^*$-algebra $A$, we denote by $A^1$ the unit boll of $A$, $A_{\rm sa}$ the set of self adjoint elements of $A$, $A_+$ the positive cone of $A$, and $U(A)$ the unitary group of $A$. Set $x\varphi(y):= \varphi(yx)$, $\varphi x(y):= \varphi(xy)$, $\|x\|_{\varphi}:= \varphi(x^*x)^{1/2}$, $\|x\|_{\varphi}^{\sharp}$ $:=((\varphi(x^*x)+\varphi(xx^*))/2)^{1/2}$, $[x,y]=xy-yx$ for $\varphi\in A^*$ and $x$, $y\in A$.

\section{Central sequences}\label{Sec2}
The following lemma is a generalization of Lemma 3.7 in \cite{CSSWW}. To prove this  we start from the fact on the strong$*$ topology. 
Let $\cM$ be a von Neumann algebra and  $x_n\in \cM^1$, $n\in\N$. If  $\|[x_n, \varphi]\|\rightarrow 0$, for any $\varphi \in \cM_*$, $x_n$, $n\in\N$ is called a {\it central sequence}. 
Let $A$ be a separable strong$*$ dense $C^*$-subalgebra of $\cM$ and $\varphi$  a faithful normal state of $\cM$. We recall that if $x_n\in \cM^1$, $n\in\N$ is a central sequence and $y_n \in \cM^1$ satisfies $\|x_n - y_n\|_{\varphi}^{\sharp} \rightarrow 0$ then $y_n \in \cM ^1$ is also a central sequence. This follows from the estimation below, (the same argument is appeared in \cite{Tak},)
\begin{eqnarray} 
\|[y_n, \psi]\| 
&\leq& \|(y_n - x_n) \psi\| + \|[x_n, \psi]\| + \|\psi(x_n - y_n)\|  \nonumber \\
&\leq& \|x_n-y_n \|_{\psi}  +\|[x_n, \psi]\| + \| ( x_n -y_n )^*\|_{\psi} \nonumber \\
&\leq& 2\sqrt{2}\|x_n-y_n\|_{\psi}^{\sharp} + \| [x_n, \psi] \| \rightarrow 0,\nonumber
\end{eqnarray}
for any $\psi\in \cM_*$.

\begin{lemma}\label{Lem0}
Let $\cM$, $\varphi$, and $A$ be as the above. Suppose that $A$ is  unital and  nuclear. Then the following holds.
\begin{description}
\item{(i)} For any central sequence $H_n \in \cM_+^1(\in \cM_{\rm sa}^1)$, $n\in\N$ there exist $h_n \in A_+^1( \in A_{\rm sa}^1)$, $n\in\N$ such that $\|h_n- H_n\|_{\varphi} \rightarrow 0$ and $\|[h_n, a]\| \rightarrow 0$ for any $a\in A$. 
\item{(ii)} For any central sequence $U_n \in U(\cM)$, $n\in\N$ there exist $u_n \in U(A)$, $n\in\N$ such that $\|u_n- U_n\|_{\varphi}^{\sharp} \rightarrow 0$ and $\|[u_n, a]\| \rightarrow 0$ for any $a\in A$. 
\end{description}
\end{lemma}
\bprf
Let $F_m$, $m\in \N$ be  finite subsets of $A^1$ such that $\overline{\bigcup F_m}^{\|\cdot\|}$ $=A^1$, and let $\varepsilon_m >0$, $m\in\N$ be a decreasing sequence such that $\varepsilon_m \searrow 0$. Because of Haagerup's theorem in \cite{Ha} we know that $A$ is amenable. And by using the approximate diagonal defined by Johnson \cite{Joh} we can obtain finite subsets $G_m$, $m\in \N$ of $A^1$ such that 

\[ \sum_{g\in G_m} g^*g =1,\quad \| [\sum_{g\in G_m} g^*ag, f]\| \leq \varepsilon_m, \]
for any $a\in A^1$ and $f\in F_m$. 

{\noindent{\it Proof of (i).\ }}
 Since $H_n \in \cM_+^1$, $n\in\N$ is a central sequence, we have a slow increasing sequence $m_n \in \N$, $n\in\N$ such that $m_n\nearrow \infty$, $|G_{m_n}|
\cdot\| [ H_n, \varphi] \|^{1/2} \rightarrow 0$, and $\displaystyle \sum_{g\in G_{m_n}} \|[ H_n, g\varphi ] \|^{1/2} \rightarrow 0$, $n\rightarrow \infty$. By  Kaplansky's density theorem, we obtain $h_n'\in A_+^1$, $n\in\N$ such that $\displaystyle \sum_{g\in G_{m_n}} \| h_n' -H_n \|_{g\varphi g^*} \leq \varepsilon_n$, $n\in\N$. Define 
\[ h_n:= \sum_{g\in G_{m_n}} g^* h_n' g\quad \in A_+^1.\]
By $\overline{\bigcup F_m}^{\|\cdot\|}=A^1$ and  the second condition of $G_m$, it follows that $\|[h_n, a]\| \rightarrow 0$ for any $a\in A^1$. Since $\sum g^*g =1$, we conclude that 
\begin{eqnarray}
\|h_n -H_n \|_{\varphi} &=& \| \sum_{g\in G_{m_n}} g^*(h_n'-H_n )g + g^*[H_n, g] \|_{\varphi} \nonumber \\
&\leq& \sum \| h_n' -H_n \|_{g\varphi g^*} + \| [H_n, g ] \|_{\varphi} \nonumber\\
&\leq& \varepsilon_n + \sum ([H_n,g]\varphi ( [H_n, g]^*))^{1/2} \nonumber \\
&\leq& \varepsilon_n + \sqrt{2} \sum (\|[H_n,g]\varphi \|)^{1/2} \nonumber \\ 
&\leq& \varepsilon_n + \sqrt{2} \sum (\|[H_n, g\varphi]\| + \|g[\varphi, H_n]\|)^{1/2} \nonumber \\
&\leq& \varepsilon_n + \sqrt{2} \sum \|[H_n, g\varphi]\|^{1/2} + \| [\varphi, H_n]\|^{1/2} \rightarrow 0. \nonumber
\end{eqnarray}
Replacing $\cM_+$ and $A_+$ with $\cM_{\rm sa}$ and $A_{\rm sa}$ in the above argument, we can verify the same. 

{\noindent{\it Proof of (ii)\ }}

Let $\log$ be the standard branch defined on the complement of the negative real axis, and $H_n:= \frac{1}{\pi \sqrt{-1}} \log (U_n)$ $\in \cM_{\rm sa}^1$, $n\in\N$. We shall show that $H_n$, $n\in\N$ is also a central sequence in the following argument.

Let  $\omega$ be a free ultrafilter on $\N$, $\mu$ be the measure on $\T$ defined by $\int_{\T} f d\mu$ $= \lim_{n\rightarrow\omega}\varphi (f(U_n))$ for any $f\in C(\T)$. Set $\|f\|_2 := (\int |f|^2 d\mu )^{1/2}$ for $f\in C(\T)$. For $\varepsilon >0$, set a polynomial $p\in C(\T)$ such that $\| p - \log \|_2 < \varepsilon /2$. Since  \begin{eqnarray} 
\|[\log(U_n), \varphi]\| &\leq& \|(\log - p) (U_n) \varphi \| + \|[p(U_n), \varphi]\| + \|\varphi (p-\log)(U_n)\| \nonumber\\
&\leq& \|[p(U_n), \varphi]\| + 2 \varphi(| \log - p |^2 (U_n))^{1/2}, \nonumber
\end{eqnarray}
 it follows that 
$\displaystyle \lim_{n\rightarrow \omega} \|[\log(U_n), \varphi]\| \leq \varepsilon$, 
which implies 
\[\lim_{n\rightarrow \omega} \|[\log(U_n), \varphi]\| =0 .\]
Then we can obtain a subsequence $n_m\in \N$ such that $\displaystyle \lim_{m\rightarrow \infty }\|[\log(U_{n_m}) , \varphi ]\| =0$. Therefore we have $\displaystyle \lim_{n\rightarrow \infty} \|[\log (U_n), \varphi]\|=0$.

Applying ${\it (i)}$ to $H_n\in \cM_{\rm sa}^1$, $n\in\N$, we obtain $h_n \in A_{\rm sa}^1$, $n\in \N$ satisfying the conditions in ${\it (i)}$. Set $\delta_n:=\max \{ \|[\varphi, h_n] \|, \| [\varphi, H_n]\|\}$. Thus we have $\delta_n \rightarrow 0$, and for any $k\in \N$ we have  
\begin{eqnarray}
\|h_n^k - H_n^k \|_{\varphi} &\leq& \|(h_n - H_n) h_n^{k-1}\|_{\varphi} + \| H_n (h_n^{k-1} - H_n^{k-1} )\|_{\varphi} \nonumber \\
&\leq& (\varphi ((h_n- H_n) h_n^{2(k-1)} (h_n- H_n)) \nonumber \\
&+& [\varphi, h_n^{k-1}(h_n- H_n)]((h_n-H_n)h_n^{k-1}))^{1/2}\nonumber \\
&+& \|h_n^{k-1}- H_n^{k-1}\|_{\varphi}\nonumber\\
&\leq& \|h_n-H_n\|_{\varphi} + \sqrt{2}\|[\varphi, h_n^{k-1} (h_n-H_n) ]\|^{1/2}+\|h_n^{k-1} -H_n^{k-1}\|_{\varphi}\nonumber\\
&\leq&\|h_n-H_n\|_{\varphi} + 2(k\delta_n)^{1/2} + \| h_n^{k-1} - H_n^{k-1}\|_{\varphi}\nonumber\\
&\leq&k\|h_n-H_n\|_{\varphi} + 2\delta_n^{1/2} \sum_{l=1}^k l^{1/2} \nonumber \\
&\leq&2^k\|h_n-H_n\|_{\varphi} + \delta_n^{1/2}k(k+1)\nonumber \\
&\leq& 2^k\|h_n-H_n\|_{\varphi} + 4^k\delta_n \nonumber 
\nonumber.
\end{eqnarray}

Define $u_n:= \exp(\pi\sqrt{-1}h_n)$, $n\in \N$. Hence it follows that $\| [u_n, a]\| \rightarrow 0$ for any $a\in A^1$ and 
\begin{eqnarray}
\|u_n - U_n \|_{\varphi}^{\sharp} &=& \|\exp (\pi\sqrt{-1} h_n) - \exp(\pi\sqrt{-1} H_n ) \|_{\varphi}^{\sharp} \nonumber \\
& \leq & \sum_{k=1}^{\infty} \frac{1}{k!} \pi^k \| h_n^k -H_n^k \|_{\varphi} \nonumber \\
& \leq &  (e^{2\pi} -1)\|h_n - H_n \|_{\varphi}+(e^{4\pi}-1)\delta_n^{1/2} \rightarrow 0.  \nonumber
\end{eqnarray}
\eprf

By ${\it (i)}$ Lemma \ref{Lem0} we obtain the following theorem which is a variation of Theorem 4.5 \cite{Kis} for positive elements and a generalization of Theorem 1.2 \cite{Sat2} for nuclear C$^*$-algebras. 

\begin{definition}
Let $A$ be a unital separable simple C$^*$-algebra with a unique tracial state $\tau$. We denote by $A_{\infty}$ the central sequence algebra $A'\cap \l^{\infty}(A)/c_0(A)$, which is defined by the operator norm topology. An automorphism $\alpha$ of $A$ has {\it the weak Rohlin property } if for any $k\in\N$ there exists $(f_n)_n \in {A_{\infty}}_+^1$ such that 
\[ (\alpha^j(f_n) f_n)_n =0,\quad j=1,2,...,k-1,\quad \lim_n \tau (1-\sum_{j=0}^{k-1} \alpha^j (f_n))=0.\]

We let  $\pi_{\tau} $ be the GNS-representation associated with $\tau$ and
\[ \WInn (A) = \{ \alpha \in \Aut (A) ; \pi_{\tau} \circ \alpha = \Ad V\circ \pi_{\tau},\ \exists V \in U(\pi_{\tau} (A)'')\}.\]
\end{definition}

\begin{theorem}\label{Thm2}
Let $A$ be a unital separable simple nuclear $C^*$-algebra with a unique tracial state $\tau$ and  $\alpha$ an automorphism of $A$.
Then  $[\alpha] \in \Aut (A) /\WInn (A)$ is aperiodic if and only if $\alpha$ has the weak Rohlin property.
\end{theorem}
\bprf
By using ${\it (i)}$  Lemma \ref{Lem0}, Theorem \ref{Thm2} becomes a trivial generalization of Theorem 1.2 in \cite{Sat2}, so the proof is sketchy.

Suppose that $\alpha \in \Aut (A) $ has a central sequence $(f_n')_n \in (A_{\infty})_+^1$ which satisfies the condition in the above definition  for $2k \in \N$. Set $f_n = \sum _{j=0}^{k-1} \alpha^j (f_n').$ Thus $(f_n)_n \in (A_{\infty})_+^1$ satisfies 
$(f_n \alpha^{k}(f_n))_n =0,$ and $\tau(1_A - (f_n + \alpha ^k (f_n))) \rightarrow 0 $. Set $\beta=\alpha^k$. If $\beta\in \WInn (A)$ then we obtain $V\in \pi_{\tau}(A)''$ such that $\pi_{\tau} \circ \beta = \Ad V \circ \pi_{\tau}$. Set a representation $\rho :A\times_{\beta} \Z \rightarrow \pi_{\tau}(A)''$, by $\rho(a)=\pi_{\tau}(a)$ and $\rho(u_{\beta})=V$, where $u_{\beta}$ is the implementing unitary for $\beta$. Then  we can define a trace $\phi$ on $A\times_{\beta} \Z$ by $\phi(x)=\tau \circ \rho (x)$, $x\in A\times_{\beta} \Z$. Since $V\in\pi_{\tau}(A)''$ we can obtain $a_n\in U(A)$, $n\in\N$ such that $\pi_{\tau}(a_n)\rightarrow V^*$ strongly, thus we have $\rho (a_nu_{\beta})\rightarrow 1$. However the partition of unities $f_n$, $\beta(f_n)$ are central sequences and transitive by $\beta$ then we have $\phi(au_{\beta})=0$ for any $a\in A$, (see Proposition 5.5 in \cite{Sat}). This contradicts $\phi (a_n u_{\beta})\rightarrow 1$.

Suppose that $[\alpha] \in \Aut(A) / \WInn (A) $ is aperiodic. Let $\balpha$ be the weak extension of $\pi_{\tau} \circ \alpha \circ \pi_{\tau}^{-1} $ to an automorphism of $\pi_{\tau} (A) ''$. 

By the classification theory for aperiodic automorphisms on the injective type $I\hspace{-.1em}I_1$ factor due to Connes \cite{Con}, there exists a central sequence $(E_m) \in (\pi_{\tau} (A) '')_{\infty}$, $m\in\N$ of projections, in the sense of the strong topology, such that 
\[E_m + \balpha(E_m)+ \cdots + \balpha^{k-1}(E_m) \rightarrow  1\ ({\rm strongly }).\]

By ${\it (i)}$ Lemma \ref{Lem0} we obtain a central sequence $(f_m')_m \in A_{\infty}$, in the operator norm  sense, such that $f_m'\in A_+^1$ and $E_m -\pi_{\tau} (f_m') \rightarrow  0 \ ({\rm strongly})$. Then we have that $\displaystyle \tau (1 - \sum_{j=0}^{k-1} \overline{ \alpha}^j \circ \pi_{\tau} (f_m')) \rightarrow 0, $ $m\rightarrow \infty$. 

By the continuous function calculus in the proof of Theorem 1.2 \cite{Sat2}, we have a orthogonality in the operator norm sense, i.e., we can take a subsequence $m_n$, $n\in\N$ and a central sequence $(f_n)_n\in A_{\infty}$ such that 
\begin{eqnarray}
 (f_n)_n (\alpha^j ( f_n ))_n &=& 0 , \quad j=1,2,...,k-1, \nonumber\\
 \|f_n- f_{m_n}'\|_{\tau} &\rightarrow& 0. \nonumber
\end{eqnarray}
Since $\displaystyle \tau (1 - \sum_{j=0}^{k-1} \overline{ \alpha}^j \circ \pi_{\tau} (f_{m_n}')) \rightarrow 0, $ $(n\rightarrow \infty)$, we also have $\displaystyle \tau (1 - \sum_{j=0}^{k-1} \alpha^j (f_n)) \rightarrow 0 $. 
\eprf

\section{Certain cocycle conjugacy result}\label{Sec3}
\begin{theorem}\label{Thm1}
Let  $\cM$ be a McDuff factor von Neumann algebra and $A$ a separable nuclear weakly dense C$^*$-subalgebra of $\cM$.
Let $G$ be a countable discrete amenable group and $\alpha$, $\beta$  centrally free actions of $G$ on $\cM$. 
Suppose that $\alpha_g \circ \beta_g^{-1} \in \oInn(\cM)$ for any $g\in G$.
Then there exist $U_g\in U(\cM)$, $g\in G$, and approximately inner automorphism $\sigma$ of $A$, in the operator norm sense, such that 
\[ \Ad U_g \circ \alpha_g = \overline{\sigma} \circ \beta_g \circ \overline{\sigma}^{-1}, \quad {\rm for\ any\ } g\in G.\]
\end{theorem}

Combining the Ocneanu's 1-cohomology vanishing theorem (7.2 in \cite{Ocn} or Lemma 6 in \cite{Mas}) with (ii) Lemma \ref{Lem0}, we obtain the following.

\begin{lemma}\label{Lem1}
Suppose that $G$, $\cM$, $A$, and $\alpha$ satisfy the conditions in Theorem \ref{Thm1}. Then for any $\varepsilon >0$ and finite subsets $F\subset G$, $\Phi \subset \cM_*$, and $E\subset A^1$, there exist $\delta >0$ and finite subset $\Psi \subset \cM_*$ satisfying the following conditions: for any $\alpha$-cocycle $U_g$ $\in U(\cM)$, $g\in G$,  with $\| [U_g, \psi]\| <\delta$, $g\in F$, $\psi
\in \Psi$, we obtain $v\in U(A)$ such that 
\[ \| U_g\alpha_g(v)^* v - 1 \|_{\varphi}^{\sharp} < \varphi, \quad \|[v, \varphi]\|< \varepsilon, \]
for any $g\in F$ and $\varphi \in \Phi$, and 
\[ \|[v, a]\| < \varepsilon, \quad a\in E. \]

Furthermore, if $\Phi$ and $E$ are empty then it may be assumed that $\Psi$ is empty.
\end{lemma}

The following lemma was appeared in Corollary 5 of \cite{Mas} which was proved by using the Ocneanu's 2-cohomology vanishing theorem in \cite{Ocn}.

\begin{lemma}\label{Lem2}
Let $\alpha$ and $\beta$ satisfy the conditions in Theorem \ref{Thm1}. Then for any $\varepsilon>0$, finite subsets $\Phi \subset \cM_*$ and $F\subset G$, there exists an $\alpha$-cocycle $U_g\in U(\cM)$, $g\in G$ such that 
\[\|\Ad U_g \circ \alpha_g (\varphi) -  \beta_g(\varphi) \| < \varepsilon, \]  
for $g\in F$ and $\varphi\in \Phi$.
\end{lemma}
\noindent{\it Proof of Theorem \ref{Thm1}.\ }

Let $A_n$, $n\in \Z_+$ be an increasing sequence of finite subsets whose union is dense in $A^1$, on the operator norm topology, and $A_0=\{1_{A}\}$. Set a faithful normal state $\varphi_0$ of $\cM$. Let $\Phi_n$, $n\in \Z_+$ be an increasing sequence of finite subsets whose union is dense in $M_*$, on the operator norm topology, and $\Phi_0 =\{ \varphi_0 \}$. 
Let $G_n$, $n\in \Z_+$ be an increasing sequence of finite subsets whose union is $G$ and $G_0= \{ 1_G \}$. 

We shall inductively construct $k_n \in \Z_+$, $U_n^{(n)} \in U(\cM)$ for $g\in G$, $v_n \in  U(A)$, $\delta_n > 0$, and finite subsets $\Psi_n$ of $\cM_*$ for $n\in \N$ satisfying the following conditions: \\
Setting $k_{-2}=k_{-1}=k_0 =0$, $U_g^{(n)} =1$ for $g\in G$, $v_0=1$, $\delta_0 =1$, $\Psi_0 =\{ \varphi_0 \}$, $\alpha^{(-1)} = \alpha$, $\beta^{(0)} =\beta$, $\oW_g^{(-1)}=\oW_g^{(0)}=1$, 
\begin{description}
\item{} $\alpha_g^{(2n+1)} = \Ad U_g^{(2n+1)} \circ \alpha_g^{(2n-1)}, \quad$ 
 $\beta_g^{(2n+2)} = \Ad U_g^{(2n+2)} \circ \beta_g^{(2n)}$, \ for \ $ g\in G$,
\item{} $\sigma_n^{(i)} =\Ad v_{2n+i}v_{2n-2+i} \cdots v_i$, \  for \ $i=0,1$,
\item{} $W_g^{(2n+1)} = U_g^{(2n+1)} \alpha_g^{(2n-1)} (v_{2n+1})^* v_{2n+1} \in \cM$, \ for \ $g\in G$, 
\item{} $W_g^{(2n+2)} = U_g^{(2n+2)} \beta_g^{(2n)} (v_{2n+2})^* v_{2n+2} \in \cM$, \ for \ $g\in G$, \ and
\item{} $\oW_g^{(2n+1+i)} = W_g^{(2n+1+i)} \Ad v_{2n+1+i}^* (\oW_g^{(2n-1+i)})\in \cM$, \ for \ $i=0,1$, $g\in G$, 
\end{description} 
the conditions are given by 

\begin{description}
\item{(1.$2n+1$)} $k_{2n+1} \geq k_{2n}$,$\quad$ $A_{k_{2n+1}} \supset_{2^{-(n+2)}} \{v_{2n}\} \cup \sigma_n^{(0)} (A_{k_{2n}})$, 
\item{(1.$2n+2$)} $k_{2n+2} \geq k_{2n}$,$\quad$ $A_{k_{2n+2}} \supset_{2^{-(n+3)}} \{v_{2n+1}\} \cup \sigma_n^{(1)} (A_{k_{2n+1}})$, 
\item{(2.2n+1)} $U_g^{(2n+1)}$, $g\in G$ is an $\alpha^{(2n-1)}$-cocycle, \\
$\| \Ad U_g^{(2n+1)} \circ \alpha_g^{(2n-1)} (\varphi) - \beta_g^{(2n)} (\varphi) \| < 2^{-1} \delta_{2n+1}$, \\
for $g\in G_{2n+1}$, $\varphi \in \bigcup_{g\in G_{2n+1}} {\alpha_g^{(2n-1)}}^{-1} (\Phi_{k_{2n+1}})\cup {\beta_g^{(2n)}}^{-1}(\Phi_{k_{2n+1}})$, 
\item{(2.$2n+2$)} $U_g^{(2n+2)}$, $g\in G$ is a $\beta^{(2n)}$-cocycle, \\
$\| \Ad U_g^{(2n+2)} \circ \beta_g^{(2n)} (\varphi) - \alpha_g^{(2n+1)} (\varphi) \| < 2^{-1} \delta_{2n+2}$, \\
for $g\in G_{2n+2}$, $\varphi \in \bigcup_{g\in G_{2n+2}} {\beta_g^{(2n)}}^{-1} (\Phi_{k_{2n+2}})\cup {\alpha_g^{(2n+1)}}^{-1}(\Phi_{k_{2n+2}})$,
\item{(3.$2n+1$)} $\| U_g^{(2n+1)} \alpha_g^{(2n-1)} (v_{2n+1})^* v_{2n+1} -1 \|_{\varphi}^{\sharp} < 2^{-(n+1)}$, \\
for $g\in G_{2n-1}$, $\varphi\in \Psi_{2n-1}$, \\
$\|[v_{2n+1}, \varphi] \| < 4^{-(n+2)}$ for $\varphi\in \Psi_{2n-1}$, \\
$\|[v_{2n+1}, a] \| < 2^{-(n+2)}$ for $a \in \sigma_{n-1}^{(1)} (A_{k_{2n-1}}) \cup A_{k_{2n-1}}$, 
\item{(3.$2n+2$)} $\| U_g^{(2n+2)} \beta_g^{(2n)} (v_{2n+2})^* v_{2n+2} -1 \|_{\varphi}^{\sharp} < 2^{-(n+2)}$, \\
for $g\in G_{2n}$, $\varphi\in \Psi_{2n}$, \\
$\|[v_{2n+1}, \varphi] \| < 4^{-(n+3)}$ for $\varphi\in \Psi_{2n}$, \\
$\|[v_{2n+2}, a] \| < 2^{-(n+3)}$ for $a \in \sigma_{n-1}^{(0)} (A_{k_{2n-2}}) \cup A_{k_{2n-2}}$, 
\item{(4.$2n+1$)} $\delta_{2n+1} \leq 2^{-1} \delta_{2n}$, and \\ 
if $\beta^{(2n)}$-cocycle $U_g\in U(\cM)$, $g\in G$, satisfies that $\| [U_g, \varphi]\| < \delta_{2n+1}$ for $g\in G_{2n+1}$, $\varphi\in \Phi_{k_{2n+1}}$, then there exists $v\in U(A)$ such that 
\[\|U_g - v^* \beta_g^{(2n)} (v) \|_{\varphi}^{\sharp} < 2^{-(n+2)}\quad {\rm for\ } g\in G_{2n+1}, \varphi\in \Psi_{2n}, \]
\[\|[v, \varphi]\| < 2^{-(n+2)}\quad {\rm for\ } \varphi\in \Psi_{2n}, \]
\[\|[v, a]\| < 2^{-(n+2)} \quad {\rm for\ } a\in A_{k_{2n}}\cup\sigma_n^{(0)} (A_{k_{2n}}), \]
\item{(4.$2n+2$)} $\delta_{2n+2} \leq 2^{-1} \delta_{2n+1}$, and \\ 
if $\alpha^{(2n+1)}$-cocycle $U_g\in U(\cM)$, $g\in G$, satisfies that $\| [U_g, \varphi]\| < \delta_{2n+2}$ for $g\in G_{2n+2}$, $\varphi\in \Phi_{k_{2n+2}}$, then there exists $v\in U(A)$ such that 
\[\|U_g - v^* \alpha_g^{(2n+1)} (v) \|_{\varphi}^{\sharp} < 2^{-(n+3)}\quad {\rm for\ } g\in G_{2n+2}, \varphi\in \Psi_{2n+1}, \]
\[\|[v, \varphi]\| < 2^{-(n+3)}\quad {\rm for\ } \varphi\in \Psi_{2n+1}, \]
\[\|[v, a]\| < 2^{-(n+3)} \quad {\rm for\ } a\in A_{k_{2n+1}}\cup\sigma_n^{(1)} (A_{k_{2n+1}}), \]
\item{(5.$2n+1$)} $\Psi_{2n+1} \supset_{2^{-(n+2)}} \Phi_{2n+1} \cup \{\oW_g^{(2n+1)} \varphi_0, \ \varphi_0\oW_g^{(2n+1)} \}_{g\in G_{2n}}$, 
\item{(5.$2n+2$)} $\Psi_{2n+2} \supset_{2^{-(n+3)}} \Phi_{2n+2} \cup \{\oW_g^{(2n+2)} \varphi_0, \ \varphi_0\oW_g^{(2n+2)} \}_{g\in G_{2n+1}}$.
\end{description}

For $n=0$, set $k_1=k_2=0$, $U_g^{(1)}=U_g^{(2)}=1$ for $g\in G$, $v_1=v_2=1$, $\delta_1=2^{-1}$, $\delta_2=2^{-2}$, and $\Psi_1=\Phi_1$, $\Psi_2=\Phi_2$. Then the conditions (1.1) $\sim$ (3.2), (5.1), and (5.2) are trivially satisfied. Since $\alpha$ and $\beta$ are centrally free actions (4.1) and (4.2) follow from Lemma \ref{Lem1}. 

Assume that we have constructed $k_{2m+1}$, $k_{2m+2}$, $U_g^{(2m+1)}$, $U_g^{(2m+2)}$, $v_{2m+1}$, $v_{2m+2}$, $\delta_{2m+1}$, $\delta_{2m+2}$, $\Psi_{2m+1}$, $\Psi_{2m+2}$ for $m\leq n-1$ which satisfy (1.$2m+1$) $\sim$ (5.$2m+2$). We proceed an induction in the following. Note that $k_m$, $U_g^{(m)}$, and $\oW_g^{(m)}$ for $m\leq 2n$ are already constructed. Since $\beta^{(2n)}$ is a centrally free action of $G$, applying Lemma \ref{Lem1} to $2^{-(n+2)}$, $G_{2n+1} \subset G$, $\Phi_{k_{2n}}$, and $A_{k_{2n}} \subset A^{1}$, we obtain $k_{2n+1}$ and $\delta_{2n+1} >0$ satisfying (1.$2n+1$) and (4.$2n+1$).  By Lemma \ref{Lem2} there exists $\alpha^{(2n-1)}$-cocycle $U_g^{(2n+1)}\in U(\cM)$ satisfying (2.$2n+1$). By (2.$2n$), i.e., $\| \alpha_g^{(2n-1)} (\varphi) - \beta_g^{(2n)} (\varphi) \| < 2^{-1} \delta_{2n}$ for $g\in G_{2n}$ and $\varphi\in \bigcup_{g\in G_{2n}} {\alpha_g^{(2n-1)}}^{-1} (\Phi_{k_{2n}})$, we have 
\begin{eqnarray}
\| [U_g^{(2n+1)}, \alpha_g^{(2n-1)} (\varphi) ] \| & =& \| \Ad U_g^{(2n+1)} \circ \alpha_g^{(2n-1)} (\varphi) - \alpha_g^{(2n-1)} (\varphi) \| \nonumber \\
&<&2^{-1}\delta_{2n+1} + 2^{-1} \delta_{2n} < \delta_{2n}. \nonumber
\end{eqnarray}
By (4.$2n$) we obtain $v^{(2n+1)}\in U(A)$ satisfying (3.$2n+1$). 
Thus we can construct $\oW^{(2n+1)}$ and obtain $\Psi_{2n+1}\subset \cM_*$ satisfying (5.$2n+1$). Now we have constructed $k_{2n+1}$, $U^{(2n+1)}$, $v_{2n+1}$, $\delta_{2n+1}$, and $\Psi_{2n+1}$ satisfying the conditions for $2n+1$. By the same way we can construct $k_{2n+2}$, $U^{(2n+2)}$, $v_{2n+2}$, $\delta_{2n+2}$, and $\Psi_{2n+2}$ satisfying the conditions for $2n+2$. 

By (1.$2n+i$) and (3.$2n+i$) for $i=0,1$ we have 
\[ \| \sigma_n^{(i)} (a) - \sigma_{n-1}^{(i)} (a) \| = \| [ v_{2n+i}, \sigma_{n-1}^{(i)} (a) ]\| < 2^{-(n+i)}, \]
for any $a\in A_{k_{2n-1+i}}$. Thus $\sigma_n^{(i)} (a)$, $n\in \N$, $i=0,1$, are Cauchy sequence in the operator norm topology for any $a\in \cup_{m} A_m$. Then we can define endomorphisms $\sigma_i$, $i=0,1$, of $A$ by 
\[ \sigma_i (a) =\lim_{n\rightarrow \infty} \sigma_n^{(i)} (a) \quad {\rm (operator\ norm)}, \]
for $a\in A$. 
By (3.$2n+i$), $i=0,1$ and $\|[v_{2n+i}, a]\| < 2^{-(n+i)}$, for $a\in A_{k_{2n-1+i}}$ we can also define endomorphisms $\sigma_i^{-1}$ of $A$ by 
\[\sigma_i^{-1} (a) = \lim_{n\rightarrow \infty} {\sigma_n^{(i)}}^{-1} (a) \quad {\rm (operator\ norm)}, \] 
for $a\in A$. It is not so hard to see that $\sigma_i\circ \sigma_i^{-1} = \id_A =\sigma_i^{-1}\circ \sigma_i$, thus $\sigma_i$ are approximately inner automorphisms of $A$.

In the following, we show that $\oW_g^{(2n+1+i)}$, $n\in\N$ are Cauchy sequences with respect to $\|\cdot \|_{\varphi_0}^{\sharp}$ for any $g\in G$ and $i=0,1$. The same calculation appeared in \cite{Mas}. By the first condition of (3.$2n+1+i$), we have 
\[\|W_g^{(2n+1+i)} -1 \|_{\varphi}^{\sharp} < 2^{-(n+1+i)}\quad {\rm for\ } g\in G_{2n-1+i}, \varphi \in \Psi_{2n-1+i}.\]
Then it follows that $\|\oW_g^{(2n+1+i)} - \oW_g^{(2n-1+i)} \|_{\varphi_0}^{\sharp}$
\begin{eqnarray}
&\leq& \| (W_g^{(2n+1+i)} -1 )(\Ad v_{2n+1+i}^* (\oW_g^{(2n-1+i)}) - \oW_g^{(2n-1+i)} \|_{\varphi_0}^{\sharp} \nonumber \\
&+& \|( W_g^{(2n+1+i)} -1 ) \oW_g^{(2n-1+i)} \|_{\varphi_0}^{\sharp} \nonumber 
\\
&+& \| \Ad v_{2n+1+i}^* (\oW_g^{(2n-1+i)}) - \oW_g^{(2n-1+i)} \|_{\varphi_0}^{\sharp}\quad  {\rm for\ any\ } g\in G.\nonumber 
\end{eqnarray}
 By (5. $2n-1+i$), $\{\Ad\oW_g^{(2n-1+i)} (\varphi_0)\}_{g\in G_{2n-2+i}}$ $\subset \Psi_{2n-1+i}$, we have 
\begin{eqnarray}
& &2{\| (W_g^{(2n+1+i)} -1 ) \oW_g^{(2n-1+i)}\|_{\varphi_0}^{\sharp}}^2 \nonumber \\
& &= \| (W_g^{(2n+1+i)} -1) \oW_g^{(2n-1+i)} \|_{\varphi_0}^2 
+ \| {\oW_g^{(2n-1+i)}}^* (W_g^{(2n+1+i)} -1)^* \|_{\varphi_0}^2 \nonumber \\
& &\leq \| W_g^{(2n+1+i)} -1\|_{\Ad\oW_g^{(2n-1+i)} (\varphi_0)}^2 + \| (W_g^{(2n+1+i)} -1 )^*\|_{\varphi_0}^2\leq 4^{-(n+i)} \nonumber
\end{eqnarray}
for $g\in G_{2n-1+i}$. By (5.$2n-1+i$) and (3.$2n+1+i$) we have
\begin{align*}
 & \| \Ad v_{2n+1+i}^* (\oW_g^{(2n-1+i)} ) - \oW_g^{(2n-1+i)} \|_{\varphi_0}^2 \nonumber \\
 & \leq \| \Ad v_{2n+1+i}^* (\oW_g^{(2n-1+i)})- \oW_g^{(2n-1+i)} \| \nonumber \\
 &\cdot \| (\Ad v_{(2n+1+i)}^* (\oW_g^{(2n-1+i)})- \oW_g^{(2n-1+i)}) \varphi_0 \| \nonumber \\
 &\leq 2 (\|v_{2n+1+i}^*\oW_g^{(2n-1+i)} [v_{2n+1+i}, \varphi_0]\| + \| v_{2n+1+i}^*[\oW_g^{(2n-1+i)} \varphi_0, v_{2n+1+i}]\| ) \nonumber \\ 
 & \leq 4^{-(n+1+i)}, \nonumber\\
\intertext{and}
& \| (\Ad v_{2n+1+i}^* (\oW_g^{(2n-1+i)} ) - {\oW_g^{(2n-1+i)}})^* \|_{\varphi_0}^2< 4^{-(n+1+i)},\quad {\rm for\ } g\in G_{2n-2+i}, \nonumber 
\end{align*}
 which implies that
\[\|\Ad v_{2n+1+i}^* (\oW_g^{(2n-1+i)})- \oW_g^{(2n-1+i)} \|_{\varphi_0}^{\sharp} < 2^{-(n+1+i)}.\] 
Combining these estimations we have 
\begin{align*}
 &{2\|(W_g^{(2n+1+i)}-1)(\Ad v_{2n+1+i}^* (\oW_g^{(2n-1+i)}) - \oW_g^{(2n-1+i)}) \|_{\varphi_0}^{\sharp}}^2 \nonumber \\
 &\leq 4 (\| \Ad v_{2n+1+i}^* (\oW_g^{(2n-1+i)}) - \oW_g^{(2n-1+i)} \|_{\varphi_0}^2 + \| (W_g^{(2n+1+i)} -1)^*\|_{\varphi_0}^2)\nonumber \\
 &\leq 8({\| \Ad v_{2n+1+i}^* (\oW_g^{(2n-1+i)}) - \oW_g^{(2n-1+i)} \|_{\varphi_0}^{\sharp}}^2 + {\| W_g^{(2n+1+i)} -1 \|_{\varphi_0}^{\sharp}}^2) \nonumber \\
&\leq 4^{-(n-1+i)}\quad {\rm for\ } g\in G_{2n-2+i}.\nonumber
\end{align*} 
Hence we conclude that $\oW_g^{(2n+1+i)}$, $n\in\N$ is a Cauchy sequence with respect to $\|\cdot\|_{\varphi_0}^{\sharp}$ for any $g\in G$.
We define 
\[ W_{i,g}= \lim_{n\rightarrow\infty} \oW_g^{(2n+1+i)} \quad {\rm (strong\ *)},\quad {\rm for\ } g\in G.\]  
From the definition of $\oW_g^{(2n+1+i)}$, $W_g^{(2n+1+i)}$, and $\sigma_n^{(i)}$ $\in \Aut (\cM)$, it follows that 
\begin{eqnarray}
\alpha_g^{(2n+1)} &=& \Ad \oW_g^{(2n+1)} \circ {\sigma_n^{(1)}}^{-1} \circ\alpha_g\circ\sigma_n^{(1)}, \nonumber \\
\beta_g^{(2n)} &=& \Ad \oW_g^{(2n)} \circ {\sigma_n^{(0)}}^{-1} \circ \beta_g \circ \sigma_n^{(0)} \quad{\rm for }\ g\in G. \nonumber 
\end{eqnarray}
 Then we can see 
\[ \Ad W_{1,g} \circ {\osigma_1}^{-1} \circ \alpha_g \circ \osigma_1 = \Ad W_{0,g} \circ {\osigma_0}^{-1} \circ \beta_g \circ \osigma_0\quad {\rm on}\ \cM,\]
for any $g\in G$. Indeed, by (5.$2n$) $\Psi_{2n} \supset \Phi_{2n-1}$ and (3.2n+1), it follows that 
\begin{eqnarray}
\|\sigma_n^{(1)}(\varphi)- \sigma_{n-1}^{(1)}(\varphi) \| &=& \| \Ad v_{2n+1} (\varphi) - \varphi \| \nonumber \\
&\leq& \|[v_{2n+1}, \varphi]\| \leq 4^{-(n+1)}\quad {\rm for\ } \varphi\in \Phi_{2n-1}, \nonumber
\end{eqnarray}
which implies $\osigma_1(\varphi) = \lim_{n\rightarrow\infty} \sigma_n^{(1)} (\varphi)$ for any $\varphi\in\cM_*$ in the norm topology. Similarly, since 
\begin{eqnarray}
\| {\sigma_n^{(1)}}^{-1}(\varphi) - {\sigma_{n-1}^{(1)}}^{-1}(\varphi) \| &=&
\| {\sigma_n^{(1)}}^{-1}\circ \sigma_{n-1}^{(1)} (\varphi) - \varphi \| \nonumber\\
&\leq& \| [v_{2n+1}, \varphi]\|, \nonumber 
\end{eqnarray}
it follows that ${\osigma_1}^{-1}(\varphi) = \lim_{n\rightarrow \infty} {\sigma_n^{(1)}}^{-1} (\varphi)$ for any $\varphi \in \cM_*$ in the norm topology. Then we have 
\[ {\osigma_1}^{-1} \circ \alpha_g \circ \osigma_1(\varphi) = \lim_n {\sigma_n^{(1)}}^{-1} \circ \alpha_g \circ \sigma_n^{(1)} (\varphi) \quad {\rm for\ any\ } \varphi\in \cM_*,\]
in the norm topology. Since 
$\|\Ad \oW_{i,g}(\varphi)- \Ad \oW_g^{(2n+i)} (\varphi) \|$
\begin{eqnarray}
&\leq& \|(\oW_{i,g} - \oW_{g}^{(2n+i)})\varphi\| + \|\varphi (\oW_{i,g} - \oW_g^{(2n+i)}) \| \nonumber \\
&\leq & \| \oW_{i,g} - \oW_g^{(2n+i)}\|_{\varphi} + \|(\oW_{i,g} - \oW_g^{(2n+i)})^* \|_{\varphi} \nonumber \\
&\leq &\sqrt{2} \| \oW_{i,g} - W_g^{(2n+i)} \|_{\varphi}^{\sharp} \rightarrow 0 \quad {\rm for\ any\ } \varphi \in \cM_*, g\in G, \nonumber
\end{eqnarray}
it follows that 
\begin{align*}
& \|\Ad \oW_{1,g}\circ {\osigma_1}^{-1}\circ \alpha_g \circ \osigma_1 (\varphi) - \Ad \oW_g^{(2n+1)} \circ {\sigma_n^{(1)}}^{-1} \circ \alpha_g \circ \sigma_n^{(1)} (\varphi) \| \nonumber\\
& \leq \|{\osigma}_1^{-1} \circ \alpha_g \circ \osigma_1 (\Ad \oW_{1,g}(\varphi)) - {\sigma_{n}^{(1)}}^{-1} \circ \alpha_g \circ \sigma_{n}^{(1)} (\Ad \oW_{1,g} (\varphi)) \| \nonumber \\
& + \| \Ad \oW_{1,g}(\varphi) - \Ad \oW_g^{(2n+1)} (\varphi) \| \rightarrow 0 \quad {\rm for}\ g\in G.\nonumber
\end{align*}
By the same way we have $\Ad \oW_{0,g} \circ {\osigma_0}^{-1} \circ \beta_g\circ \osigma_0(\varphi)$ $=$ $\lim_{n\rightarrow\infty} \Ad\oW_g^{(2n)}\circ {\sigma_{2n}^{(0)}}^{-1} \circ \beta_g \circ \sigma_{2n}^{(0)} (\varphi)$. Hence we conclude 
\[ \Ad W_{1,g} \circ {\osigma_1}^{-1} \circ \alpha_g \circ \osigma_1 (\varphi) = \Ad W_{0,g} \circ {\osigma_0}^{-1}\circ \beta_g \circ \osigma_{0} (\varphi)\]
for any $g\in G$ and $\varphi \in \cM_*$ which means the same equality on $\cM$. Setting $\sigma = \sigma_1\circ \sigma_0^{-1} \in \Aut (A)$ and $U_g= \osigma_1 (W_{0,g}^*W_{1,g}) \in U(\cM)$ for $g\in G$, we have the desired conditions. 
\eprf

\begin{remark}
By the condition in Theorem \ref{Thm1}, we can modify $U_g\in U(\cM)$ as an $\alpha$-cocycle. Then, by using the 1-cohomology vanishing lemma we can make $\alpha$-cocycle $U_g$ close to $1_{\cM}$. 

When $\alpha$ and $\beta$ make $A$ invariant, i.e., $\alpha_g(A)=\beta_g(A)=A$ for any $g\in G$,  $\Ad U_g$ is automatically an weakly inner automorphism of $A$ for any $g\in G$. Then, in particular we obtain the following corollary.
\end{remark}
\begin{corollary}
Let $A$ be a unital separable simple nuclear $C^*$-algebra with a unique tracial state $\tau$ and  $\alpha$, $\beta$ automorphisms of $A$.
If $[\alpha]$, $[\beta]$ $\in \Aut (A) /\WInn (A)$ are aperiodic then they are conjugate.
\end{corollary}
{\bf Acknowledgments.}
The author would like to thank Reiji Tomatsu who suggested that the main result  Theorem \ref{Thm1} can be extended to discrete amenable groups. He is also grateful to Hiroyuki Osaka for the travel supports at Ritsumeikan University. The preliminary work of this paper was done when he visited Sophus Lie Conference Center in Nordfjordeid Norway and  East China Normal University in Shanghai China, in June 2010. He is also  grateful to people there for their warm hospitality.

\end{document}